# On Resolving Singularities

## 0. Introduction

Let $V$ be an irreducible affine algebraic variety over a field $k$ of characteristic zero, and let $(f_0, ..., f_m)$ be a sequence of elements of the coordinate ring. There is probably no elementary condition on the $f_i$ and their derivatives which determines whether the blowup of $V$ along $(f_0, ..., f_m)$ is nonsingular. The result of the current paper is that there indeed is such an elementary condition, involving the first and second derivatives of the $f_i$, provided we admit certain singular blowups, all of which can be resolved by an additional Nash blowup. This paper is the promised sequel of [1], in which the same program was carried out for individual vector fields. Indeed, the current paper generalizes the result of [1] to algebraic foliations of arbitrary codimension, and the case of codimension zero foliations corresponds to the problem of resolving the singularities of $V$.

## 1. The Nash Resolution Question

The editors have asked for an account of Nash's question and the relation with this work. Here is Milnor's account from [3], where further references may be found. Let $r = dim(V)$. Suppose $V = V_0 \subset W_0$ is an embedding in a nonsingular variety over $k$. Then $V_0$ lifts to a subvariety $V_1 \subset W_1 = Grass_r(W_0)$ of the variety of $r$-planes in the tangent bundle of $W_0$. The natural map $\pi : V_1 \to V_0$ is called the Nash blowup of $V_0$. It is the lowest blowup where $\pi^*(\Omega_{V_0/k})$/torsion is locally free. Now we can repeat the process, giving a variety $V_2 \subset W_2 = Grass_r(W_1)$ and so-on, and the question is whether eventually $V_i$ is nonsingular. There is a particular explicit sequence of ideals $R = J_0, J_1, J_2, ... \subset R$ so that $V_0 = Bl_{J_0}V, V_1 = Bl_{J_1}V, V_2 = Bl_{J_2}V, ....$ with $J_i|J_{i+1}$ for all $i$. Applying our earlier paper [2], $V_i$ is nonsingular if and only if the ideal class of $J_{i+1}$ divides some power of the ideal class of $J_i$. The present paper brings things down to earth considerably: such a divisibility of ideal classes implies that for some $N \geq r + 2$ and $i$

$$J_i^{N-r-2}J_{i+1}^{r+3} = J_i^N J_{i+2}.$$

Yet note that this identity in turn implies $J_{i+2}$ is a divisor of some power of $J_{i+1}$. Thus although $V_i$ may fail to be nonsingular, when the identity holds the *next* variety $V_{i+1}$ must be nonsingular. Thus the Nash question is equivalent to the assertion that the identity above holds for some sufficiently large $i$ and $N$.

John Atwell Moody
Maths Institute, Warwick University, Coventry CV4 7AL, England
moody@maths.warwick.ac.uk
Math Reviews Subject Classification: 14A05 (secondary 32C45)

## 2. A toy theorem

In order to explain the main theorem, Theorem 15, let's look first at what it says about an individual vector field $\delta$ on $V$ when $k = \mathbb{C}$, $V$ is an irreducible algebraic variety in $\mathbb{C}^3$, and $(f_0, ...., f_m)$ is the trivial (unit) ideal. Let $R$ be the coordinate ring of $V$ and let $K$ be the function field of $V$. Choose any nonzero vector field on $V$, corresponding to a $k$ linear derivation

$$\delta : R \to R$$

$$r \mapsto \dot{r}.$$

Let $(x, y, z)$ be the standard coordinates in $\mathbb{C}^3$ and let us define the velocity and acceleration vectors by

$$v = (\dot{x}, \dot{y}, \dot{z})$$

$$\dot{v} = (\ddot{x}, \ddot{y}, \ddot{z}).$$

Let $(v)$ be the ideal generated by the entries of $v$ and let $(v \times \dot{v})$ be the ideal generated by the entries of the cross product. Let $S \subset V$ be the vanishing locus of $(\dot{x}, \dot{y}, \dot{z})$ and let $\widetilde{V} \subset \mathbb{C}^3 \times \mathbb{P}^2$ be the closure of the graph of the familiar Gauss map

$$V - S \to \mathbb{P}^2$$

$$(x, y, z) \mapsto [\dot{x} : \dot{y} : \dot{z}].$$

Note that if $\delta$ induces a nonsingular foliation on $V$ then it also induces a nonsingular foliation on $\widetilde{V} = V$. Now, in this very special situation we have **Toy Theorem**.

1. (i) There is always an inclusion $(v)^3 \subset (v)^3 + (v \times \dot{v})$.

2. (ii) If $\delta$ induces a nonsingular foliation on $V$ then the inclusion above is an equality $(v)^3 = (v)^3 + (v \times \dot{v})$.

3. (iii) Conversely, if the inclusion is an equality then $\delta$ induces a nonsingular foliation on $\widetilde{V}$.

Note that if $V$ is a curve then $\delta$ induces a nonsingular foliation on $V$ (respectively $\widetilde{V}$) if and only if $V$ (respectively $\widetilde{V}$) is nonsingular.

Theorem 15 is an analogous theorem which works not only for $V$ itself, but for an arbitrary blowup of $V$, and for foliations of any dimension. Moreover, the ideals which occur in the statment of Theorem 15 are all ideals of the original ring $R$. The proof of the toy theorem works like this. Part (i) is obvious. Let $\pi : \widetilde{V} \to V$ be the natural map. For any sequence of elements $l = (l_0, ..., l_n)$ of $K$ we may consider the fractional ideal $(l) \subset K$, which is the set of $R$ linear combinations of the $l_i$. For any such $(l)$ we define a new fractional ideal $\mathcal{J}(l) = (l \times \dot{l}) + (l)^2(v)$ (where $(l \times \dot{l})$ is a fractional ideal with $\binom{n+1}{2}$ entries).

Some properties of $\mathcal{J}(l)$ which were proven in [1] are

1. 1. $\pi^*\mathcal{J}(v) = \mathcal{O}_{\widetilde{V}}\delta(\mathcal{O}_{\widetilde{V}})(-2E)$
2. 2. $(l)^2\mathcal{J}(m) + (m)^2\mathcal{J}(l) = \mathcal{J}(lm)$
3. 3. $\mathcal{J}(1) = (v)$ where $\pi^*$ is the operation of pulling back ideals to ideal sheaves, $lm$ is the product sequence and $E$ is the exceptional divisor associated to blowing up $(v)$.

Here is the proof of (ii). Suppose $\delta$ induces a nonsingular foliation on $V$. This means the fractional ideal $(v)$ is invertible, ie there is a sequence $(w)$ of elements of $K$ so that $(v)(w) = R$. Applying properties 2 and 3 we have

$$(w)^2(v \times \dot{v}) \subset (w)^2\mathcal{J}(v) \subset \mathcal{J}(vw) = \mathcal{J}(1) = (v).$$

Multiplying both sides by $(v)^2$ gives $(v \times \dot{v}) \subset (v)^3$ proving the desired equality of ideals $\mathcal{J}(v) = (v)^3$. Here is the proof of (iii). Suppose conversely that there is such an equality of ideals $\mathcal{J}(v) = (v)^3$. Applying $\pi^*$ and applying 1 yields

$$\mathcal{O}_{\widetilde{V}}\delta(\mathcal{O}_{\widetilde{V}})(-2E) \cong \mathcal{O}_{\widetilde{V}}(-3E)$$

so $\mathcal{O}_{\widetilde{V}}\delta(\mathcal{O}_{\widetilde{V}}) \cong \mathcal{O}_{\widetilde{V}}(-E)$ is locally free proving $\delta$ induces a nonsingular foliation on $\widetilde{V}$. The main thing which the toy theorem is meant to illustrate is that the explicit and elementary condition that the inclusion $(v)^3 \subset (v)^3 + (v \times \dot{v})$ should be an equality is intermediate in strength between saying $\delta$ induces a nonsingular foliation on $V$ and saying $\delta$ induces a nonsingular foliation on $\widetilde{V}$. Theorem 15 concerns our arbitrary irreducible affine variety $V$ over our field $k$ of characteristic zero, with coordinate ring $R$, an arbitrary ideal $I = (f_0, ..., f_m) \subset R$, and an arbitrary algebraic foliation $L$ on $V$ (of any dimension). The theorem describes an elementary condition involving first and second derivatives, which is intermediate in strength between saying $L$ lifts to a nonsingular foliation on $Bl_I(V)$ and saying $L$ lifts to a nonsingular foliation on the Gauss blowup $\widetilde{Bl_I(V)}$ of $Bl_I(V)$ along $L$. An important aspect of Theorem 15 is that the ideals considered in the statement are all ideals in the original ring $R$. Now I will give the statement of Theorem 15. Let $R$ be the coordinate ring of $V$, and let $x_1, ..., x_n$ be a sequence of $k$ algebra generators of $R$. Recall $K$ is the fraction field of $R$. An algebraic (singular) foliation on $V$ corresponds to a $K$ sub Lie algebra $L \subset \mathcal{D}er_k(K, K)$ (see section 3 below); let $\delta_1, ..., \delta_r$ be a $K$ basis of $L$. We may assume $\delta_i \in \mathcal{D}er_k(R, R)$. For any ideal $I = (f_0, ..., f_m)$ of $R$ let $\mathcal{J}(I)$ be the ideal generated by the

$$f_{u_1}f_{u_2}...f_{u_b} \cdot \det \begin{pmatrix} f_{i_1} & \delta_1 f_{i_1} & ... & \delta_r f_{i_1} \\ & ... & & \\ f_{i_a} & \delta_1 f_{i_a} & ... & \delta_r f_{i_a} \\ 0 & \delta_1 x_{j_1} & ... & \delta_r x_{j_1} \\ 0 & ... & & \\ 0 & \delta_1 x_{j_b} & ... & \delta_r x_{j_b} \end{pmatrix}$$

where $a$ and $b$ run over numbers such that $a + b = r + 1$, and where $0 \leq u_1, ..., u_b \leq m, 0 \leq i_1 < i_2 < ... < i_a \leq m$, and $1 \leq j_1 < j_2 < ... < j_b \leq n$. The ideal $\mathcal{J}(I)$ is independent of choice of generators $(f_0, ..., f_m)$. For each choice of ideal $I$ the ideal $\mathcal{J}(I)$ has the property that the Gauss blowup $\widetilde{Bl_I(V)}$ of $Bl_I(V)$ along $L$ satisfies $\widetilde{Bl_I(V)} = Bl_{I\mathcal{J}(I)}(V)$ (see Lemma 4 below). For each choice of ideal $I$ the ideal $\mathcal{J}(I\mathcal{J}(I))$ makes sense, and it is generated by certain explicit expressions involving the $f_i$ and the $x_j$ and their first and second derivatives. Theorem 15 in this situation makes the following three assertions:

**Main Theorem.**

1. (i) For any ideal $I$ of $R$ there is an inclusion $\mathcal{J}(I)^{r+2} \subset \mathcal{J}(I\mathcal{J}(I))$ of ideals of $R$.

2. (ii) If $L$ lifts to a nonsingular foliation on $Bl_I(V)$ then this inclusion becomes an equality after multiplying both sides by a suitable $N$'th power of $I$:
$$I^N \mathcal{J}(I)^{r+2} = I^N \mathcal{J}(I\mathcal{J}(I)).$$

3. (iii) Conversely if the equality in (ii) holds then $L$ does lift to a nonsingular foliation on the Gauss blowup $\widetilde{Bl_I(V)}$ of $Bl_I(V)$ along $L$.

When $L$ is the unique codimension zero foliation of $V$ then the lift of $L$ to any blowup of $V$ is just the unique codimension zero foliation of the blowup, which we may also call $L$. To say $L$ is nonsingular on a blowup is the same as saying the blowup is a nonsingular variety. Moreover, the Gauss blowup of $Bl_I(V)$ along $L$ is just the Nash blowup of $Bl_I(V)$. Assembling these facts we see that if $I$ resolves the singularities of $V$, so $Bl_I(V)$ is a nonsingular variety, then for $r = dim(V)$ and some $N$ the inclusion $I^N \mathcal{J}(I)^{r+2} \subset I^N \mathcal{J}(I\mathcal{J}(I))$ becomes an equality, and conversely when this is so the Nash blowup $\widetilde{Bl_I(V)}$ is a nonsingular variety.

3. Singular Foliations

Let $V$ be an affine irreducible variety over a field $k$. $V$ is determined by its coordinate ring, an arbitrary finite type $k$-algebra $R$ without zero divisors. Let $K$ be the fraction field of $R$. Let's say a *singular foliation* on $V$ is just any $K$ linear Lie sub algebra $L \subset \mathcal{D}er_k(K, K)$. We shall let $\widehat{L} = Hom_K(L, K)$ and we shall let $\overline{\Omega_{V/k}}$ be the image of the homomorphism

$$\Omega_{V/k} \to \widehat{L}$$

which sends a differential of the form $dx$ to the functional $(\delta \mapsto \delta(x))$. Let us record this formula to avoid any confusion:

$$\overline{\Omega_{V/k}} = Image(\Omega_{V/k} \to \widehat{L}).$$

We shall say the singular foliation is *nonsingular on $V$* if this image $R$-module is projective. The following proposition justifies this definition. **1. Proposition.**

If $L$ is nonsingular on $V$ then there are elements $f_i \in R$ generating the unit ideal such that for each $i$ each Lie ring $L \cap \mathcal{D}er_k(R[f_i^{-1}], R[f_i^{-1}])$ is free over $R[f_i^{-1}]$ with a basis $\delta_1, ..., \delta_r$ (depending on $i$), such that are elements $x_1, ..., x_r \in R[f^{-1}]$ with $\delta_i(x_j) = 1, i = j$, and $\delta_i(x_j) = 0, i \neq j$. Proof. Apply $Hom_R(-, R) \subset Hom_R(-, K)$ to the split surjection $\Omega_{V/k} \to \overline{\Omega_{V/k}}$ to obtain a pullback square showing
$$L \cap \mathcal{D}er_k(R, R) = Hom_R(\overline{\Omega_{V/k}}, R).$$

We may assume $\overline{\Omega_{V/k}}$ is free with basis of the form $dx_1, ..., dx_r$, and let $\delta_1, ..., \delta_r$ be the dual basis.

Even if $V$ is a non-affine irreducible variety with function field $K$, for $L$ a $K$-linear Lie algebra $L \subset \mathcal{D}er_k(K, K)$, we may still think of $L$ as defining a singular foliation on $V$. $L$ is said to be nonsingular just if it is nonsingular on each affine part, and the proposition above still holds for each affine part of $V$. If $W$ is any variety birationally equivalent to $V$, and if $L \subset \mathcal{D}er_k(K, K)$ defines a singular foliation on $V$, then note $L$ also defines a singular foliation on $W$ because $V$ and $W$ share the same function field. In particular if $W$ is a blowup of $V$, we may consider the question of whether the singular foliation $L$ on $V$ becomes nonsingular on $W$. In what follows we'll restrict to the case that $V$ is affine with coordinate ring $R$, furnished with a singular foliation $L$, and we will study ideals $I \subset R$ for the question of whether there exists an ideal $I$ such that $L$ becomes nonsingular on $\widetilde{V} = Bl_I(V)$.

## 4. The Module $M_\gamma$

The previous section leads to the question of how to compute $\overline{\Omega_{\widetilde{V}/k}}$ in terms of $\overline{\Omega_{V/k}}$ and the ideal $I$. In these terms there is a natural answer: there is a class depending on $I$

$$\gamma \in Ext^1_R(I, I \otimes \Omega_{V/k})$$

which has a natural image, which we'll also call $\gamma$, in each $Ext^1_R(I, I\overline{\Omega_{V/k}})$, and this class defines a certain $R$-module $M_\gamma$ with a structure map $p: M_\gamma \to I$. (The referee has written that the class $\gamma$ is the well-known 'Atiyah class' (see section 4 of [3]) and therefore the module $M_\gamma$ is equal to an appropriate quotient of the module of principal parts of $I$. In other words $M_\gamma$ is an appropriate quotient of the module of global 1-jets of sections of $I$ viewed as a coherent sheaf.) In Proposition 2 below we'll show that there is an exact sequence of sheaves on $\widetilde{V} = Bl_I(V)$

$$0 \to \overline{\Omega_{\widetilde{V}/k}} \to \overline{\pi^*(M_\gamma)}(E) \xrightarrow{\overline{\pi^*p(E)}} \mathcal{O}_{\widetilde{V}} \to 0 \qquad (1)$$

where $E$ is the exceptional divisor. Here the overline in the middle term refers to reduction modulo torsion. $\pi: \widetilde{V} \to V$ is the structure map of the blowups. The map $\overline{\pi^*p(E)}$ is the result of pulling back via $\pi$ the map $p: M_\gamma \to I$, twisting by by $E$, and reducing modulo torsion. It is perhaps a good time to explain the convention which we'll use throughout this paper. An overline on a module or coherent sheaf will always denote the torsion-free quotient of that module or sheaf, *except* that for any variety $V$ over $k$ the symbol $\overline{\Omega_{V/k}}$ will denote the natural image of of $\Omega_{V/k}$ in $\widehat{L}$, which we'll call the *reduced* differentials, and, more generally, an overline over $\wedge^r \Omega_{V/k}$ will denote the natural image in $\wedge^r \widehat{L}$. Let us also display this more general formula to avoid any confusion:

$$\overline{\wedge^r \Omega_{V/k}} = Image(\wedge^r \Omega_{V/k} \to \wedge^r \widehat{L}).$$

Finally, when $V$ is not affine we have a similar definition for twisted sheaves which will play a role in section 5. When $L$ is the whole of $\mathcal{D}er_k(K, K)$ these notions coincide, the torsion free quotient being the image in $\widehat{L}$. The class $\gamma$ is defined as follows: take a resolution of $I$

$$R^a \to R^b \to R^c \xrightarrow{\epsilon} I \to 0$$

and say the middle map is

$$g_i \mapsto \sum_j a_{ij} h_j$$

where $g_i$ and $h_j$ are respective basis vectors. Then a cocycle $z: R^b \to I\overline{\Omega_{V/k}}$ representing $\gamma$ sends $g_i$ to $-\sum_j \epsilon(h_j)\overline{d}a_{ij}$, where $\overline{d}$ is the natural derivation $R \to \overline{\Omega_{V/k}}$. Leibniz rule implies that this map satisfies the cocycle condition. We can build the extension module $M_\gamma$ using the double complex

$$\begin{array}{ccccccccc} R^a & \to & R^b & \to & R^c & \to & I & \to & 0 \\ \downarrow & & \downarrow z & & \downarrow e & & \downarrow & & \\ 0 & \to & I\overline{\Omega_{V/k}} & \subset & \overline{\Omega_{V/k}} & \to & 0 & & \end{array}$$

where the vertical map $e$ is given

$$e(h_j) = \overline{d}\epsilon h_j.$$

It is easy to see that the diagram commutes. A copy of $M_\gamma$ occurs as the submodule of $I \oplus \overline{\Omega_{V/k}}$ generated by the image of $I\overline{\Omega_{V/k}} + R^c$. This in turn is equal to the submodule of $I \oplus \widehat{L}$ generated by the $f \oplus \overline{d}f$ for $f \in I$. The projection $I \oplus \widehat{L} \to I$ induces a surjection $M_\gamma \to I$ and the kernel is exactly $I\overline{\Omega_{V/k}}$. One caution is that if $(f_0, ..., f_n)$ is a sequence of generators of $I$, it does not automatically follow that the rows $(f_i, \overline{d}f_i)$ generate $M_\gamma$. It is only the case that these rows together with a system of generators of $I\overline{\Omega_{V/k}}$ suffice.

**2. Proposition.** The sequence (1) is exact; i.e., the kernel of $\overline{\pi^* p(E)}$ is the image in $\widehat{L}$ of the sheaf of differentials of $\widetilde{V}$. Proof. Choose once and for all a $K$ basis $\delta_1, ..., \delta_r$ of $L$ and so we have $\widehat{L} \cong K^r$ by which each element $v$ is sent to $(v(\delta_1), ..., v(\delta_r))$. The extension module $M_\gamma$ is then isomorphic to the submodule of $I \oplus K^r$ generated by all rows $(f, \delta_1 f, ..., \delta_r f)$ for $f \in I$. Explicitly, if $I$ is generated by $f_0, ..., f_n$, we may view $M_\gamma$ as the module of $I \oplus K^r$ generated by rows of the following two types:

$$\begin{pmatrix} f_i & \delta_1(f_i) & ... & \delta_r(f_i) \\ 0 & f_j \delta_1(x_c) & ... & f_j \delta_r(x_c) \end{pmatrix}.$$

The sheaf $\overline{\pi^* M_\gamma}(E)$ can be constructed chart by chart. The 0'th coordinate chart of $\widetilde{V}$ is $U = Spec(\widetilde{R})$ for $\widetilde{R} = R[f_1/f_0, ..., f_n/f_0]$ and the module $\overline{\pi^* M_\gamma}(E)(U)$ over this ring can be explicitly constructed within $K \oplus K^r$ by multiplying each row above by $f_0^{-1}$ and considering the $\widetilde{R}$ module the new rows generate. Some typical rows which result are listed below

$$\begin{pmatrix} 1 & f_0^{-1} \delta_1(f_0) & ... & f_0^{-1} \delta_r(f_0) \\ f_i/f_0 & f_0^{-1} \delta_1(f_i) & ... & f_0^{-1} \delta_r(f_i) \\ 0 & \delta_1(x_c) & ... & \delta_r(x_c) \\ 0 & f_j/f_0 \delta_1(x_c) & ... & f_j/f_0 \delta_r(x_c) \end{pmatrix}.$$

Subtracting $f_i/f_0$ times the first row from the second, and $f_j/f_0$ times the third row from the fourth yields

$$\begin{pmatrix} 1 & f_0^{-1}\delta_1(f_0) & ... & f_0^{-1}\delta_r(f_0) \\ 0 & f_0^{-1}\delta_1(f_i) - f_i/f_0^2 \delta_1(f_0) & ... & f_0^{-1}\delta_r(f_i) - f_i/f_0^2 \delta_r(f_0) \\ 0 & \delta_1(x_c) & ... & \delta_r(x_c) \\ 0 & 0 & ... & 0 \end{pmatrix}$$

The kernel of the projection on $\widetilde{R} \oplus 0$ is generated by rows such as the second and third displayed, and using the rule for differentiating a quotient we see these are just the images of the $k$ algebra generators $x_c$ and $f_i/f_0$ of $\widetilde{R}$. The $\widetilde{R}$ module they span is the image of the natural map

$$\Omega_{\widetilde{R}/k} \to \widehat{L}$$

which is $\overline{\Omega_{\widetilde{V}/k}}(U)$ as claimed. The same considerations apply to each other coordinate chart, and this proves the kernel of the projection $\overline{\pi^* M_\gamma}(E) \to \mathcal{O}_{\widetilde{V}}$ is $\overline{\Omega_{\widetilde{V}/k}}$. It follows from the proposition that the reduced differentials of $\widetilde{V}$ are locally free if and only if the pullback of $M_\gamma$ modulo torsion is locally free, so the question of resolving the singular foliation $L$ on $V$ comes down to finding an $I$ such that the associated extension module $M_\gamma$ pulls back to a locally free sheaf modulo torsion.

There is a lowest blowup $\pi$ which makes $\overline{\pi^* M_\gamma}$ projective, namely the blowup of the rank one torsion free module $\overline{\wedge^{r+1} M_\gamma}$. We may make a fractional ideal $\mathcal{J}(I)$ isomorphic to this module; namely the fractional ideal generated by the determinants of all possible matrices

$$\begin{pmatrix} f_0 & \delta_1 f_0 & ... & \delta_r f_0 \\ & & ... & \\ f_r & \delta_1 f_r & ... & \delta_r f_r \end{pmatrix}$$

where $\delta_i$ are our fixed basis of $L$ and $f_0, ..., f_r$ ranges over all possible lists of $r+1$ elements of $I$. If the $\delta_i$ are chosen to belong to $L \cap \mathcal{D}er_k(R, R)$ then $\mathcal{J}(I)$ will be an ordinary ideal instead of a fractional ideal, but this is an unimportant limitation because we will want to apply the operator $\mathcal{J}$ to fractional ideals anyway. Note that because of proposition 2 (the exactness of the sequence (1))

we have **3. Corollary.**

$$\overline{\pi^* \mathcal{J}(I)} \cong \overline{\pi^* \wedge^{r+1} M_\gamma} \cong \overline{\wedge^r (\Omega_{\widetilde{V}/k}(-E))} \otimes \mathcal{O}_{\widetilde{V}}(-E) \cong \overline{\wedge^r \Omega_{\widetilde{V}/k}}(-E - rE).$$

Because twisting does not affect blowing up, the blowup of $\widetilde{V} = Bl_I(V)$ along $\overline{\pi^* \mathcal{J}(I)}$ is the same as the blowup of $\widetilde{V}$ along $\overline{\wedge^r \Omega_{\widetilde{V}/k}}$. This is in turn isomorphic as a variety over $V$ to the Gauss blowup of $Bl_I(V)$ along $L$ – let's call this $\widetilde{Bl_I(V)}$. Thus we have the following isomorphisms of varieties over $V$:

$$\widetilde{Bl_I(V)} = Bl_{\overline{\wedge^r \Omega_{Bl_I(V)/k}}} Bl_I(V)$$

$$\cong Bl_{\overline{\wedge^r \Omega_{Bl_I(V)/k}}(-E-rE)} Bl_I(V)$$

$$\cong Bl_{\overline{\pi^* \mathcal{J}(I)}} Bl_I(V)$$

$$\cong Bl_{I\mathcal{J}(I)} V.$$

Let's record this as a little lemma. **4. Lemma.** The Gauss blowup of $Bl_I(V)$ along $L$ is isomorphic as a variety over $V$ to $Bl_{I\mathcal{J}(I)}(V)$. **5. Proposition.** If $I$ is generated by $f_0, ..., f_m$ then $\mathcal{J}(I)$ is generated by the elements

$$f_{u_1} f_{u_2} ... f_{u_b} \cdot \det \begin{pmatrix} f_{i_1} & \delta_1 f_{i_1} & ... & \delta_r f_{i_1} \\ & ... & & \\ f_{i_a} & \delta_1 f_{i_a} & ... & \delta_r f_{i_a} \\ 0 & \delta_1 x_{j_1} & ... & \delta_r x_{j_1} \\ 0 & & ... & \\ 0 & \delta_1 x_{j_b} & ... & \delta_r x_{j_b} \end{pmatrix}$$

where $a$ and $b$ run over numbers such that $a + b = r + 1$, and where $0 \le u_1, ..., u_b \le m, 0 \le i_1 < i_2 < ... < i_a \le m$, and $1 \le j_1 < j_2 < ... < j_b \le n$. *Proof.*

In the first display in the proof of proposition 2 we saw that the image of $M_\gamma$ is generated by the rows

$$( f_i \quad \delta_1(f_i) \quad ... \quad \delta_r(f_i) )$$

and the rows

$$( 0 \quad f_j \delta_1(x_c) \quad ... \quad f_j \delta_r(x_c) ).$$

The size $r + 1$ square matrix above is obtained by choosing $r + 1$ such rows in all possible ways and taking determinants. Therefore the determinants generate the corresponding image of $\wedge^{r+1} M_\gamma$ in $K$.

The problem of resolving the singular foliation $L$ on $V$ comes down to finding an ideal $I \subset R$ such that $\overline{\wedge^r \Omega_{Bl_I(V)/k}}(-E - rE) \cong \overline{\pi^* \mathcal{J}(I)}$ is locally free. This happens if and only if $Bl_I(V)$ dominates $Bl_{\mathcal{J}(I)}(V)$. Luckily we have the following theorem: **6. Theorem** (Ill. J. Math, to appear [2]): Let $I, J \subset R$ be ideals. Then $Bl_I(V)$ dominates $Bl_J(V)$ if and only if there is a number $\alpha$ and a fractional ideal $S$ such that

$$JS = I^\alpha.$$

Although the theorem is stated for ideals, it follows for fractional ideals. Thus

**7. Corollary.** An ideal $I \subset R$ has the property that $L$ lifts to a nonsingular foliation on $\widetilde{V} = Bl_I(V)$ if and only if $\mathcal{J}(I)$ is a divisor of a power of $I$ as a fractional ideal, ie if and only if there is a fractional ideal $S$ of $R$ and a number $\alpha$ such that

$$S\mathcal{J}(I) = I^\alpha.$$

*Proof.* This is just a matter of assembling data already proven. By definition $L$ lifts to a nonsingular foliation on $\widetilde{V} = Bl_I(V)$ if and only if $\overline{\Omega_{\widetilde{V}/k}}$ is locally free. This happens if and only if $\overline{\wedge^r \Omega_{\widetilde{V}/k}}$ is locally free. And we have from Corollary 3 that $\overline{\wedge^r \Omega_{\widetilde{V}/k}} \cong \overline{\pi^* \mathcal{J}(I)}(E + rE)$. This is locally free if and only if the blowup $Bl_I(V)$ dominates $Bl_{\mathcal{J}(I)}(V)$, and by Theorem 6 this happens if and only if there is a fractional ideal $S$ of $R$ and a number $\alpha$ so $\mathcal{J}(I)S = I^\alpha$. If we take the case $r = 1$ as a guide [1], we should not expect to have a formula which will simply give us the generating sequence $(f_1, ..., f_m)$ for an ideal which will resolve the singularities of $V$. However, we may hope to have an elementary condition on the $f_i$ and their derivatives which will tell us whether the ideal or an associated ideal will resolve.

## 5. Calculation of the Reduced Differentials of the Blowup

The idea of this section, which is independent of the rest of the paper, is to describe the differentials $\Omega_{\widetilde{V}/k}$ of a blowup $\pi : \widetilde{V} \to V$ or rather the image $\overline{\Omega_{\widetilde{V}/k}} \subset \widehat{L}$. When $L = \mathcal{D}er_k(K,K)$ this is just the torsion free quotient of the differentials of $\widetilde{V}$. In the following section we'll return to the problem of determining nonsingularity of the foliation lifted to the blowup solely in terms of the generators of the ideal downstairs. During this section, though, we'll allow ourselves to work with sheaves upstairs in the blowup. Throughout this section $V$ will be affine and irreducible over a field $k$, with coordinate ring $R$. The first step is to notice for any ideal $I \subset R$ letting $\Omega = \Omega_{V/k}$ that canonical derivation $d : R \to \Omega_{V/k}$ defines a *homomorphism* of modules

$$h : I \to \Omega/I\Omega$$

because for $f \in I$ and $r \in R$ we have $d(rf) = rdf + fdr$ and the second term is in $I\Omega$. Letting $\widetilde{V} = Bl_I(V)$ there are natural inclusions of sheaves

$$\overline{\pi^*\Omega}(-E) \subset \overline{\Omega_{\widetilde{V}/k}}(-E) \subset \overline{\pi^*\Omega}$$

which define $\overline{\Omega_{\widetilde{V}/k}}(-E)$ as the inverse image in $\overline{\pi^*\Omega}$ of a certain subsheaf $\mathcal{L} \subset \overline{\pi^*\Omega}/\overline{\pi^*\Omega}(-E)$. We shall describe the sheaf $\mathcal{L}$. **8. Theorem.** The desired subsheaf $\mathcal{L}$ is the image of the composite

$$(\pi^*(I))_{tors} \subset \pi^*(I) \xrightarrow{\pi^*h} \pi^*(\Omega/I\Omega) \to \overline{\pi^*\Omega}/\overline{\pi^*\Omega}(-E).$$

Proof. Recall there is a certain $R$-module $M_\gamma$ which I defined earlier, and which fits into a pullback diagram

$$\begin{array}{ccccccccc}
 & & 0 & & 0 & & & & \\
 & & \downarrow & & \downarrow & & & & \\
0 & \to & \overline{I\Omega} & \to & \overline{\Omega} & \to & \overline{\Omega}/\overline{I\Omega} & \to & 0 \\
 & & \downarrow & & \downarrow & & \downarrow & & \\
0 & \to & M_\gamma & \to & I \oplus \overline{\Omega} & \to & \frac{I \oplus \overline{\Omega}}{M_\gamma} & \to & 0 \\
 & & \downarrow & & \downarrow & & & & \\
 & & I & = & I & & & & \\
 & & \downarrow & & \downarrow & & & & \\
 & & 0 & & 0 & & & &
\end{array}$$

and defines an isomorphism $\overline{\Omega}/\overline{I\Omega} \to \frac{I \oplus \overline{\Omega}}{M_\gamma}$. This diagram of course does not remain a pullback after pulling back by $\pi$, but the important property of $M_\gamma$ is that we do obtain a similar diagram by applying $\pi^*$ just to the maps $M_\gamma \to I$ and $I \oplus \overline{\Omega} \to I$, reducing mod torsion and taking kernels.

$$\begin{array}{ccccccccc}
& & 0 & & 0 & & & & \\
& & \downarrow & & \downarrow & & & & \\
0 & \to & \overline{\Omega_{\widetilde{V}/k}}(-E) & \to & \overline{\pi^*\Omega} & \to & \overline{\pi^*\Omega}/\overline{\Omega_{\widetilde{V}/k}}(-E) & \to & 0 \\
& & \downarrow & & \downarrow & & \downarrow & & \\
0 & \to & \overline{\pi^*M_\gamma} & \to & \mathcal{O}_{\widetilde{V}}(-E) \oplus \overline{\pi^*\Omega} & \to & (\mathcal{O}_{\widetilde{V}}(-E) \oplus \overline{\pi^*\Omega})/(\overline{\pi^*M_\gamma}) & \to & 0 \\
& & \downarrow & & \downarrow & & & & \\
& & \mathcal{O}_{\widetilde{V}}(-E) & = & \mathcal{O}_{\widetilde{V}}(-E) & & & & \\
& & \downarrow & & \downarrow & & & & \\
& & 0 & & 0 & & & &
\end{array}$$

This gives us an isomorphism $\overline{\pi^*\Omega}/\overline{\Omega_{\widetilde{V}/k}}(-E) \to (\mathcal{O}_{\widetilde{V}}(-E) \oplus \overline{\pi^*\Omega})/(\overline{\pi^*M_\gamma})$. We now compare the middle row of the second diagram with the result of pulling back the middle row of the first diagram, but making substitutions according to the two isomorphisms we have so far discovered. We obtain the diagram with exact rows and surjective vertical maps

$$\begin{array}{ccccccc}
& \pi^*M_\gamma & \to & \pi^*I \oplus \pi^*\Omega & \to & \pi^*(\overline{\Omega}/I\overline{\Omega}) & \to & 0 \\
& \downarrow & & \downarrow & & \downarrow & & \\
0 \to & \overline{\pi^*M_\gamma} & \to & \mathcal{O}_{\widetilde{V}}(-E) \oplus \overline{\pi^*\Omega} & \to & \overline{\pi^*\Omega}/\overline{\Omega_{\widetilde{V}/k}}(-E) & \to & 0
\end{array}.$$

The sequence of first two kernels splices to the sequence of two displayed cokernels to give the four term exact sequence

$$(\pi^*M_\gamma)_{tors} \to (\pi^*I)_{tors} \oplus (\pi^*\overline{\Omega})_{tors} \to \pi^*(\overline{\Omega}/I\overline{\Omega}) \to \overline{\pi^*\Omega}/\overline{\Omega_{\widetilde{V}/k}}(-E) \to 0.$$

Taking into account that

$$\frac{\pi^*(\overline{\Omega}/I\overline{\Omega})}{Image((\pi^*\overline{\Omega})_{tors})} \cong \overline{\pi^*\Omega}/\overline{\pi^*\Omega(-E)}$$

gives us the exact sequence

$$(\pi^*I)_{tors} \to \overline{\pi^*\Omega}/\overline{\pi^*\Omega(-E)} \to \overline{\pi^*\Omega}/\overline{\Omega_{\widetilde{V}/k}}(-E) \to 0.$$

Here the subsheaf described explicitly in the statement of the theorem is the image of the leftmost map, and the sheaf $\mathcal{L}$ which determines the differentials of $\widetilde{V}$ is the kernel of the second map. The fact that they are equal follows therefore by exactness. **9. Corollary.** We may reconstruct $\overline{\Omega_{\widetilde{V}/k}(-E)}$ as the pullback

$$\begin{array}{ccc}
\overline{\Omega_{\widetilde{V}/k}(-E)} & \to & \overline{\pi^*\Omega} \\
\downarrow & & \downarrow \\
\mathcal{L} & \subset & \overline{\pi^*\Omega}/\overline{\pi^*\Omega(-E)}
\end{array}.$$

where $\mathcal{L}$ is the sheaf explicitly described in the statement of theorem 10.

## 6. Study of $\mathcal{J}$

We shall study the effect of $\mathcal{J}$ on powers $I^n$. We first have a lemma which, although not necessary for the main result, nevertheless has led to a simplification. I wish to thank D. Rumynin for pointing out this proposition. **10. Proposition.** Let $R$ be a $k$ algebra of finite type, $k$ a field of characteristic zero, $I$ an ideal of $R$. Let $r \geq 0$. Then there is a system of generators $(f_0, ..., f_m)$ of $I$ such that $(f_0^{r+1}, ..., f_m^{r+1}) = I^{r+1}$. Proof. Start with a sequence $(f_0, ..., f_{m'})$ which generates $I$. Then $I^{r+1}$ is generated by all degree $r+1$ monomials in the $f_i$. Each such monomial $f_0^{i_0} f_1^{i_1} ... f_{m'}^{i_{m'}}$ is equal to $1/n!$ times the alternating sum of the $r+1$'st powers of the subsums of the expression

$$\underbrace{f_0 + ... + f_0}_{i_0 \ times} + .... + \underbrace{f_{m'} + ... + f_{m'}}_{i_{m'} \ times}.$$

The additional generators can be taken to be the subsums. At this stage it turns out to be a good idea to define, for any sequence of elements $(f_1, ..., f_m)$ in $I$, the "wrong" fractional ideal $\mathcal{M}(f_1, ..., f_m)$ to be generated by the determinants resulting from these generators only

$$det \begin{pmatrix} f_{i_0} & \delta_1 f_{i_0} & ... & \delta_r f_{i_0} \\ & & ... & \\ f_{i_r} & \delta_1 f_{i_r} & ... & \delta_r f_{i_r} \end{pmatrix}.$$

This is contained in $\mathcal{J}(I)$ but the inclusion is proper in general. However, if the generating sequence $(f_0, ..., f_m)$ is appropriately enlarged the inclusion becomes an equality: **11. Proposition.** If we begin with a sequence of generators $f_0, ..., f_m$ of an ideal $I$ of $R$ and extend by appending all products with a system of $k$ algebra generators of $R$, the new sequence $(f_0, ..., f_{m'})$ will have the property that $\mathcal{J}(I) = \mathcal{M}(f_1, ..., f_{m'})$. Proof. It suffices to show the rows $(f_i \ \delta_1(f_i) \ ... \ \delta_r(f_i))$ generate $M_\gamma$. We know $M_\gamma$ is generated by rows of the type above together with rows $(0 \ f_i\delta_1(x_j) \ ... \ f_i\delta_r(x_j))$, both with $i \leq m$. Rows of the first type belong to our proposed generating set. To obtain rows of the second type, choose $s$ such that $f_s = x_j f_i$. Then

$$(f_s \ \delta_1(f_s) \ ... \ \delta_r(f_s)) - x_j(f_i \ \delta_1(f_i) \ ... \ \delta_r(f_i))$$
$$= (0 \ \delta_1(x_j f_i) - x_j \ \delta_1(f_i) \ ... \ \delta_r(x_j f_i) - x_j \ \delta_r(f_i))$$
$$= (0 \ f_i\delta_1(x_j) \ ... \ f_i\delta_r(x_j))$$

as needed.

Suppose $V$ is an irreducible variety over a field $k$ of characteristic zero, $R$ is the coordinate ring of $V$, and $K$ is the function field of $V$. Let $L \subset \mathcal{D}er_k(K, K)$ be a singular algebraic foliation on $V$. Fix $\delta_1, ..., \delta_r$ a $K$-basis of $L$. **12. Theorem.**

Let $I$ and $J$ be fractional ideals of $R$.

$$I^{(r+1)}\mathcal{J}(J) \subset \mathcal{J}(IJ).$$

Proof. Choose a sequence of generators of each ideal. Since the characteristic of $k$ is zero we can include enough generators $f_i$ of $I$ so that the powers $f_i^{r+1}$ generate $I^{r+1}$. Extend both generating sequences by appending all multiples of the generators with the $k$-algebra generators $x_i$ of $R$. Call the new sequences $(f_0, ..., f_m)$ and $(g_0, ..., g_u)$. Note that the product sequence $(f_s g_t)$ contains a system of generators of $IJ$ as well as all products of these generators with the $x_i$. Therefore by Proposition 11 we may write

$$\mathcal{M}(f_0, ..., f_m) = \mathcal{J}(I)$$
$$\mathcal{M}(g_0, ..., g_u) = \mathcal{J}(J)$$
$$\mathcal{M}((f_0, ..., f_m)(g_0, ..., g_u)) = \mathcal{J}(IJ).$$

By our choice of generators of $I$ we also have by proposition 12

$$(f_0^{r+1}, ..., f_m^{r+1}) = I^{r+1}.$$

A typical generator of $\mathcal{M}((f_0, ..., f_m)(g_0, ..., g_u))$ is a determinant of a size $r+1$ matrix whose rows look like

$$(f_s g_t, f_s \delta_1 g_t + g_t \delta_1 f_s, ..., f_s \delta_r g_t + g_t \delta_r f_s)$$

for various choices of $s$ and $t$. Each column after the first of such a matrix is a sum of two columns in an obvious way. The determinant is therefore a sum of the $2^r$ determinants where we have chosen either the first or second column in each case. If we choose the same value of $s$ for each row, all these determinants vanish except one, which is the determinant of a matrix whose rows look like

$$(f_s g_t, f_s \delta_1 g_t, ..., f_s \delta_r g_t)$$

for various values of $t$. All the other matrices in the sum have one column at least which is a multiple of the first column by an element of $K$, so their determinants vanish. The determinant of the one matrix which counts is $f_s^{r+1}$ times an arbitrary generator of $\mathcal{M}(g_0, ..., g_u)$. Repeating the calculation for each value of $s$ gives

$$(f_0^{r+1}, ..., f_m^{r+1})\mathcal{M}(g_0, ..., g_u) \subset \mathcal{M}((f_0, ..., f_m)(g_0, ..., g_u)).$$

Combining facts we have

$$I^{r+1}\mathcal{J}(J) = (f_0^{r+1}, ..., f_m^{r+1})\mathcal{J}(J) = (f_0^{r+1}, ..., f_m^{r+1})\mathcal{M}(g_0, ...., g_u)$$
$$\subset \mathcal{M}((f_0, ..., f_m)(g_0, ..., g_u)) = \mathcal{J}(IJ).$$

Applying this result plus induction we have

$$I^{(r+1)(N-1)}\mathcal{J}(I) \subset \mathcal{J}(I^N).$$

We have thus bounded $\mathcal{J}(I^N)$ from below. What is remarkable is that we may bound it from above, and the bounds will be equal, so we will have calculated $\mathcal{J}(I^N)$. Both bounding arguments will have used the fact that the characteristic of $k$ is zero, for two different reasons.

**13. Lemma.** Let $I$ be any fractional ideal of $R$ for $R$ as above. Then $\mathcal{J}(I^N) \subset I^{(N-1)(r+1)}\mathcal{J}(I)$. Proof. Let $f_0, ..., f_m$ be a generating sequence of $I$ chosen by Proposition 8 so $\mathcal{M}(f_0, ..., f_m) = \mathcal{J}(I)$. Note that the sequence of degree N monomials in the $f_i$ becomes extended at the same time in the appropriate way to satisfy the hypothesis of Proposition 8 so $\mathcal{J}(I^N)$ is equal to $\mathcal{M}$ applied to the sequence of degree $N$ monomials in the $f_i$. The latter is generated by the determinants of certain matrices. Let's now look at the case $N = 3$, the general case being similar. Each row of the typical matrix looks like $(f_if_jf_k, \ \delta_1(f_if_jf_k), \ ... \ , \delta_r(f_if_jf_k))$. Expanding out using the Leibniz rule one obtains a sum of three rows, namely

$$f_if_j(f_k/3, \ \delta_1(f_k), \ ... \ , \delta_r(f_k)),$$

$$f_if_k(f_j/3, \ \delta_1(f_j), \ .... \ , \delta_r(f_j)),$$

and

$$f_jf_k(f_i/3, \ \delta_1(f_i), \ ... \ , \delta_r(f_i)).$$

Because of the multilinearity of the determinant, our expression for the determinant is $1/3$ times a sum of degree $2(r+1)$ monomials in the $f_i$ times determinants of the matrices which come into the definition of $\mathcal{M}(f_0, ..., f_m)$. Thus each term is an element of $I^{2(r+1)}\mathcal{M}(f_0, ..., f_m) = I^{2(r+1)}\mathcal{J}(I)$ as needed. The proof clearly generalizes to arbitrary N. By some miracle, the lemmas above are precise converses of each other, so we get an equality of ideals, at least when the characteristic of $k$ is zero. **14. Theorem.** Suppose $char(k) = 0$. Let $I$ be a fractional of $R$ and let $r$ be the dimension of $L$ over $K$. Then

$$\mathcal{J}(I^N) = I^{(N-1)(r+1)}\mathcal{J}(I).$$

We will use Theorem 12 and Theorem 14 in the proof of Theorem 15. There again it will happen that separate arguments will furnish upper and lower bounds for an ideal, and these will match exactly.

# 7. The Main Theorem

Let $R$ be an integral domain which is a $k$ algebra for $k$ a field of characteristic zero. Let $V = Spec(R)$. Let $L \subset \mathcal{D}er_k(K,K)$ be a $K$-linear sub Lie algebra, and let $\delta_1, ..., \delta_r$ be a basis of $L$. For a fractional ideal $J$ of $R$ recall that $\mathcal{J}(J)$ is the fractional ideal of $R$ generated by the determinants

$$det \begin{pmatrix} f_0 & \delta_1 f_0 & ... & \delta_r f_0 \\ & & ... & \\ f_r & \delta_1 f_r & ... & \delta_r f_r \end{pmatrix}$$

for $f_0, ..., f_r \in J$. An explicit finite list of generators of $\mathcal{J}(J)$ is given in proposition 5 if $f_0, ..., f_m$ generate $J$. We can arrange that $\mathcal{J}(J)$ is an ordinary ideal instead of a fractional ideal if we bother to choose the $\delta_i$ to lie in $\mathcal{D}er_k(R,R)$, but this is an unimportant distinction. Indeed we'll end up working with fractional ideals during the proof of Theorem 15 anyway. By our definition $L$ describes a nonsingular foliation on the blowup $\widetilde{V} = Bl_J(V)$ if and only if the sheaf

$$\overline{\Omega_{\widetilde{V}/k}} = Image(\Omega_{\widetilde{V}/k} \xrightarrow{\overline{d}} \widehat{L})$$

is locally free on $\widetilde{V}$, where $\overline{d}$ is the map sending a generating section $dx$ to the function

$$L \to K$$
$$\delta \mapsto \delta(x).$$

Moreover we have proven in Corollary 7 that $L$ does define a nonsingular foliation on $\widetilde{V}$ if and only if there is a fractional ideal $S$ for $R$ and a number $\alpha$ such that

$$S\mathcal{J}(J) = J^\alpha.$$

The most important case of this is when $L = \mathcal{D}er_k(K,K)$ in which case $L$ defines a nonsingular foliation on $\widetilde{V}$ if and only if $\widetilde{V}$ is nonsingular. **15. Theorem.**

1. (i) There is always an inclusion $\mathcal{J}(J)^{r+2} \subset \mathcal{J}(J\mathcal{J}(J))$.

2. (ii) If $L$ lifts to a nonsingular foliation on the blowup $\widetilde{V} = Bl_J(V)$ then there is an $N$ such that the inclusion becomes an equality after multiplying both sides by $J^N$; i.e.,

$$J^N \mathcal{J}(J)^{r+2} = J^N \mathcal{J}(J\mathcal{J}(J)). \quad (*)$$

3. (iii) Suppose conversely that $J$ is any ideal such that the inclusion in (i) becomes an equality as in (ii). Then letting $I = J\mathcal{J}(J)$ we have that $L$ lifts to a nonsingular foliation on the blowup $Bl_I(V)$, which is the same as the Gauss blowup $\widetilde{Bl_J(V)}$ of the variety $Bl_J(V)$ along $L$.

Proof. For part (i) we have
$$\mathcal{J}(J)^{r+2} = \mathcal{J}(J)^{r+1}\mathcal{J}(J).$$
By Theorem 12 we have
$$\mathcal{J}(J)^{r+1}\mathcal{J}(J) \subset \mathcal{J}(J\mathcal{J}(J)).$$
Combining gives the result Now for the proof of (ii). Suppose $L$ lifts to a nonsingular foliation on $Bl_J(V)$. By Corollary 7 this means there is a fractional ideal $S$ and a number $\alpha$ so that
$$\mathcal{J}(J)S = J^\alpha.$$
Now we have by Theorem 12
$$S^{r+1}\mathcal{J}(J\mathcal{J}(J)) \subset \mathcal{J}(SJ\mathcal{J}(J)).$$
Combining, we see
$$S^{r+1}\mathcal{J}(J\mathcal{J}(J)) \subset \mathcal{J}(J^{\alpha+1}).$$
Using Theorem 14,
$$\mathcal{J}(J^{\alpha+1}) = J^{(r+1)\alpha}\mathcal{J}(J).$$
Combining the last two formulas and multiplying through by $\mathcal{J}(J)^{r+1}$ gives
$$(\mathcal{J}(J)S)^{r+1}\mathcal{J}(J\mathcal{J}(J)) \subset J^{(r+1)\alpha}\mathcal{J}(J)^{r+2}.$$
Again applying the result of Corollary 7 we see that $(\mathcal{J}(J)S) = J^\alpha$. Substituting this in the left side of the displayed equation gives
$$J^{(r+1)\alpha}\mathcal{J}(J\mathcal{J}(J)) \subset J^{(r+1)\alpha}\mathcal{J}(J)^{r+2}.$$
Setting $N = (r+1)\alpha$ gives one desired inclusion of ideals; the opposite inclusion is part (i), which has already been proven, multiplied by $J^N$. The combination gives the equality of ideals
$$J^N\mathcal{J}(J\mathcal{J}(J)) = J^N\mathcal{J}(J)^{r+2}.$$
This is a second time in the paper that two unrelated arguments give upper and lower bounds for an ideal, and the bounds match exactly. Now for the proof of part (iii). Suppose the equality above holds. Let $I = J\mathcal{J}(J)$. We may assume $N \geq r+2$ and let $\beta = N - r - 2$. Multiplying both sides of the formula by $\mathcal{J}(J)^\beta$ we have
$$\mathcal{J}(J)^\beta J^N \mathcal{J}(I) = J^N \mathcal{J}(J)^{\beta+r+2} = (J\mathcal{J}(J))^N = I^N.$$
Letting $S = \mathcal{J}(J)^\beta J^N$ we have
$$S\mathcal{J}(I) = I^N.$$
By corollary 7 this proves $L$ lifts to a nonsingular foliation on $Bl_I(V)$. Finally, identify $Bl_I(V)$ with the Gauss blowup of $Bl_J(V)$ along $L$ by Lemma 4.

## 8. Connection with the Nash Resolution Question

Let us connect the theorem of the previous section, in the case of the unique codimension zero foliation $L$, with the Nash question. Recall that $V$ is affine irreducible over $k$ a field of characteristic zero, $K$ the function field of $V$, $R$ its coordinate ring. For $L = \mathcal{D}er_k(K,K)$, we let $\delta_1, ..., \delta_r$ be a $K$ basis of $L$, which we can assume lie in $\mathcal{D}er_k(R,R)$. We defined for each ideal $J$ of $R$ a new ideal $\mathcal{J}(J)$ of $R$ with the property that $\mathcal{J}(J)$ is a fractional ideal divisor of a power of $J$ if and only if $L$ lifts to a nonsingular foliation on $\widetilde{V} = Bl_J(V)$. Moreover by Proposition 1, since $L$ is the unique codimension zero foliation, this happens if and only if $Bl_I(V)$ is nonsingular. Recall by Lemma 4 the blowup of the product $J\mathcal{J}(J)$ is the same as the result of blowing up $J$ to get $\widetilde{V}$ and then blowing up the highest exterior power of the reduced differentials of $\widetilde{V}$. Theorem 15, the main theorem of the previous section, in this situation says that when $\widetilde{V}$ is nonsingular, so the second blowup is an isomorphism, then the inclusion $J^N \mathcal{J}(J)^{r+2} \subset J^N \mathcal{J}(J\mathcal{J}(J))$ becomes an equality for some $N$, and that when this equality does hold then the result of blowing up the highest exterior power of the reduced differentials of $\widetilde{V}$ has a nonsingular foliation defined by $L$. One can consider a chain of ideals

$$J_0 = R$$
$$J_1 = \mathcal{J}(R)$$
$$J_2 = \mathcal{J}(R)\mathcal{J}(\mathcal{J}(R))$$
$$J_3 = \mathcal{J}(R)\mathcal{J}(\mathcal{J}(R))\mathcal{J}(\mathcal{J}(R)\mathcal{J}(\mathcal{J}(R)))$$
$$...$$
$$J_{i+1} = J_i \mathcal{J}(J_i).$$

The result of blowing up $J_i$ is the same as starting with $V$ and sequentially blowing up the highest exterior power of the reduced differentials, to obtain a sequence of varieties $V_i \to V_{i-1} \to ... \to V_0 = V$ Thus, blowing up the ideal $J_i$ accomplishes in one step what could otherwise be done by $i$ steps of blowing up the highest exterior power of the reduced differentials:

$$V_i = Bl_{J_i} V.$$

The chain of blowups stops (with all higher blowups being isomorphisms) if and only if some $J_i$ resolves the singularities of $V$. After this, the ideal classes of the higher ideals $J_{i+1}, J_{i+2}, ...$, which are clearly multiples of $J_i$, are also divisors of a power of $J_i$. In this context, the result of the previous section tells you how to check when you have successfully resolved $V$. It says that one needs only check that the inclusion

$$J_i^N \mathcal{J}(J_i)^{r+2} \subset J_i^N \mathcal{J}(J_{i+1})$$

is an equality for some $N$. When $J_i$ resolves, this condition holds, and when the condition holds $J_{i+1}$ resolves. Since $L = \mathcal{D}er_k(K, K)$ the sequence of blowups above is just the sequence of 'Nash' blowups, and the Nash question asks whether they eventually resolve the singularities of $V$. Therefore, we have a completely explicit reformulation of the Nash conjecture: **16. Theorem.** The Nash ques-tion holds in the affirmative for $V$ if and only if the inclusion above becomes an equality for some sufficiently large $i$ and $N$. To obtain the formulation in section

1, note when the inclusion above is an equality it remains so when both sides are multiplied by $J_{i+1}$ and apply the basic definitions.